\documentclass[12pt]{amsart}
\usepackage{amssymb,amscd}
\usepackage{graphicx, color}

\newtheorem*{Whitney towers}{Theorem~\ref{Whitney towers}}
\newtheorem*{h-towers}{Theorems ~\ref{half} \& \ref{$(n)$-solvable}}

\newtheorem*{surgery curves}{Theorem~\ref{surgery curves}}
\newtheorem*{cg=0}{Theorem~\ref{vanish}}

\theoremstyle{definition}

\numberwithin{equation}{section}
\numberwithin{figure}{section}

\newcommand{\bb}{\bigbreak}

\newcommand{\np}{\newpage}

\newcommand{\R}{\mathbb{R}}

\def\yen{{\setbox0=\hbox{Y}Y\kern-.97\wd0\vbox{hrule height.lex width.98%
\wd0\kern.33ex\hrule height.lex width.98\wd0\kern.45ex}}}

\def\np{\newpage}

\begin{document}  
\pagestyle{plain}

\title{Make your Boy surface}
\author{Eiji Ogasa}

\date{}


\begin{abstract} 
This is an introductory article about the Boy surface. 
Boy found  1901 that $\R P^2$ can be immersed into $\R^3$, and 
published it.  
(The image of) the immersion is called the Boy surface after Boy's discovery. 

We have created a way to construct the Boy surface 
by using a pair of scissors, a piece of paper, and a strip of scotch tape. 
In this article we introduce the way. 

Furthermore, we make a movie to show the paper-craft actually, 
and put it  in a website.    
One can find the website by typing in 
the author's name or 
the title of this article     
in the search engine. 
\end{abstract} 
\maketitle

\section{Introduction}\label{Introduction} 
This is an introductory article about the Boy surface. 
We draw the Boy surface in the following page.

Boy discovered 1901 that $\R P^2$ can be immersed into $\R^3$, and published it 
in \cite{Boy}.  
Boy is the name of a mathematician in Germany, May 4, 1879 -- September 6, 1914. 
The immersion of $\R P^2$ into $\R^3$ which Boy found  
is called the Boy surface after his discovery.
Note: The Boy surface is the name of (an image of) an immersion not that of a manifold.

See \cite{MilnorStasheff} for mathematical terms: $\R P^2$, immersions, etc. 
P121 of  \cite{MilnorStasheff} quotes \cite{Boy}. 

\cite{O0} is the author's introductory book. 

\bigbreak
We  have created a way to construct the Boy surface 
by using a pair of scissors, a piece of paper, and a strip of scotch tape. 
In this article we introduce the way.

\bigbreak 
In \cite{Giller} Giller used the Boy surface. 
In \cite{O3} the author cited  \cite{Giller}, used the Boy surface, and proved a theorem. 
In \cite{O4} the author cited  \cite{Giller} again. 
It is his motivation to write this article.

\bigbreak 
The author made a movie to show  our paper-craft actually. 
He put it in Youtube. 
The website of the movie is connected with the author's website. 
One can find these websites by typing in 
the author's name, `Eiji Ogasa', 
or 
the title of this article, `Make your Boy surface',   
in the search engine. 

\bigbreak 
The author translated this article into Japanese (\cite{Oni}).
One can find it by typing in 
the author's name, `Eiji Ogasa', or the Japanese title.

\np

\vskip-40mm
\includegraphics[width=13cm]{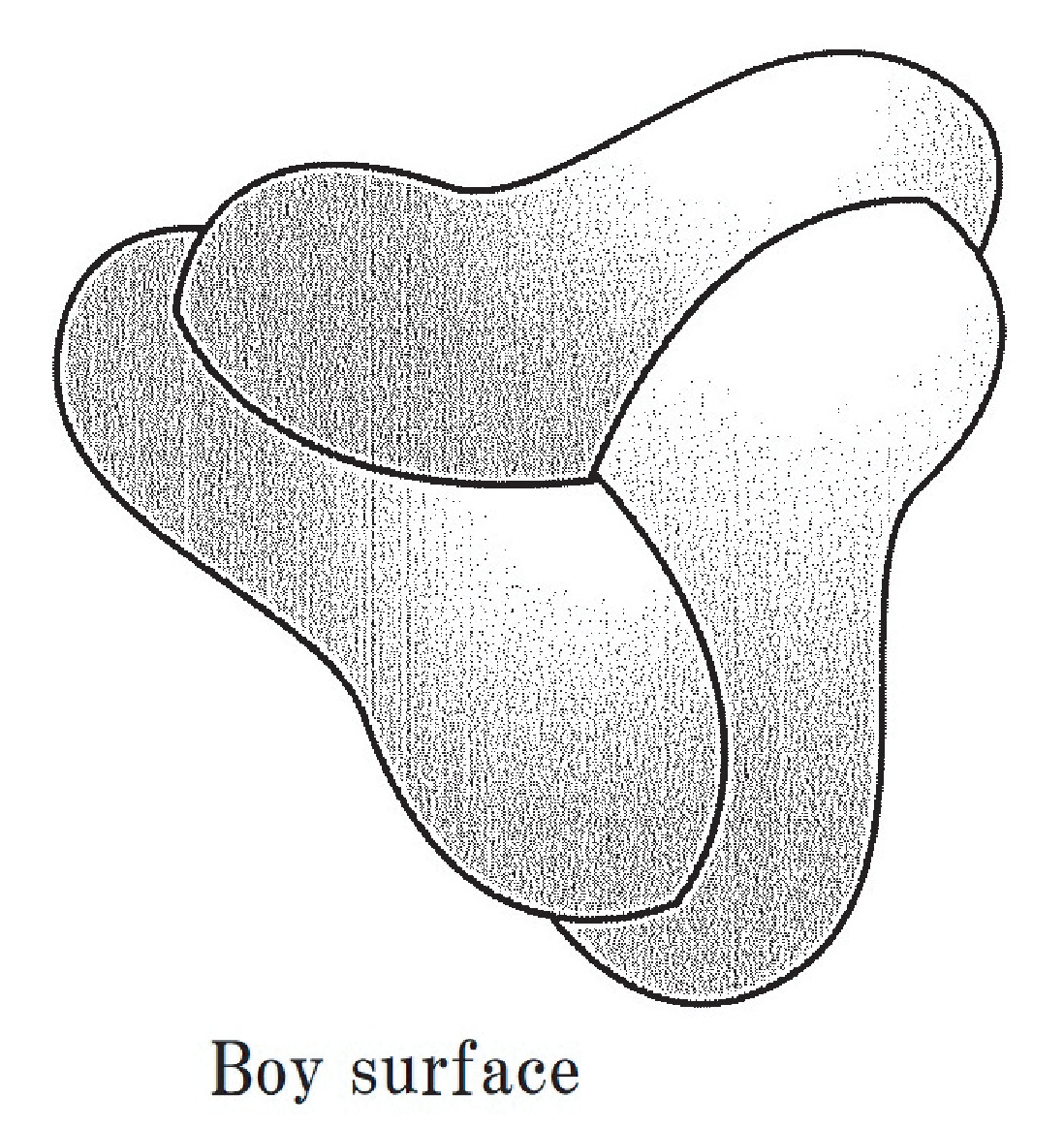} 

\np\section{Paper-Craft}\label{Paper-Craft}
See  Figure I, Figure II, and Figure III in the following two pages. 
Make three copies of Figure I, a copy of Figure II, and three copies of Figure III.

\bigbreak 
\noindent
Note:  Make the copies so that 
the length of the edge of each of the unit squares in Figure I is 
half of that  in Figure II, III.  
If it might be difficult to take such a copy of Figure II (resp. III), 
then we recommend the following way: 
Take a copy of Figure I at first.
After that, make Figure II, III on a paper by using a scale and a pencil.

\np

\vskip-40mm
\includegraphics[width=13cm]{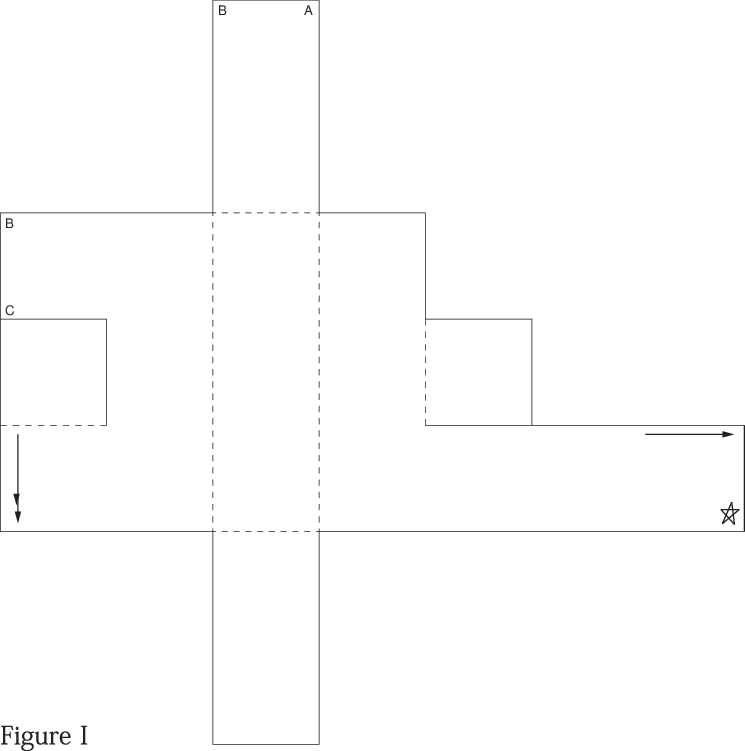}

\np

\vskip-40mm
\includegraphics[width=13cm]{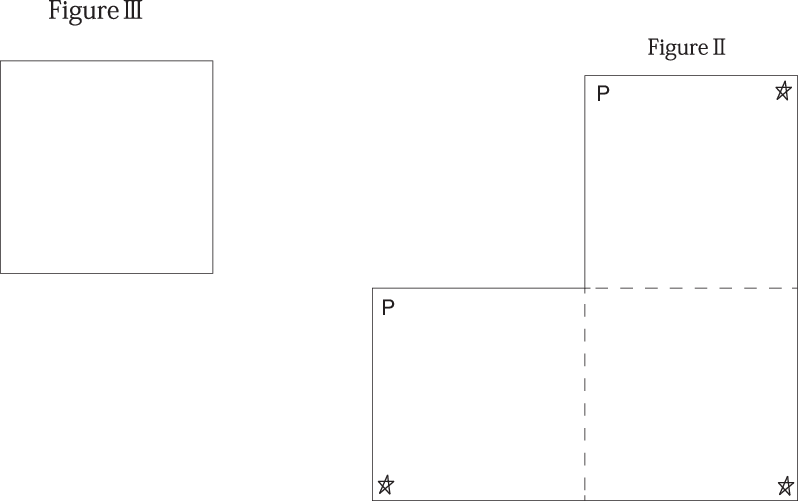}

\np

Make Figure IV from the three copies of Figure III. What we obtain is called the piece IV. 
Note that we must cut the three copies of  Figure III a few times. If necessarily, we divide 
one of the three copies into a few pieces once and  
after that, we make the piece IV from them by using a strip of scotch tape. 

\vskip1cm
The piece IV is represented as follows if we take the $xyz$-axises .  It is a union of 

$\{(x,y,z)  \vert -1\leqq x\leqq1,\quad -1\leqq y\leqq1, \quad z=0\}$

and 

$\{(x,y,z)  \vert -1\leqq y\leqq1,\quad -1\leqq z\leqq1, \quad x=0\}$

and 

$\{(x,y,z)  \vert -1\leqq z\leqq1,\quad -1\leqq x\leqq1, \quad y=0\}$.

\vskip1cm
Take the following points in the piece IV as shown in Figure V. 

$A(-1, 0,0), \quad B(-1,0,1), \quad C(0,0,1)$

$A'(0,-1,0), \quad B'(1,-1,0), \quad C'(1,0,0)$

$A''(0,0,-1), \quad B''(0,1,-1), \quad C''(0,1,0)$

We will use these points soon.

\np

\vskip-80mm
\hskip-10mm\includegraphics[width=10cm]{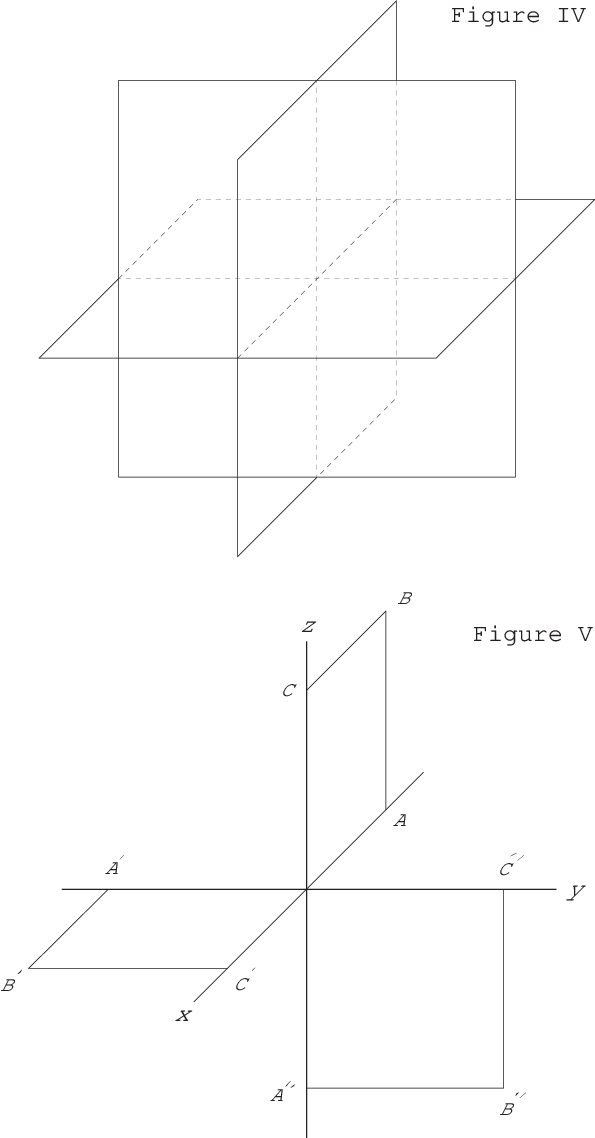}  

\np

Cut each copy of Figure I along the solid lines. What we obtain after cutting is called 
the piece I. 

Fold the piece I along the dotted line so that we see the dotted line inside, 
and make `the angle made by the paper at the dotted line' $90^\circ$.  

Use a strip of scotch tape and attach the edges which meet. 
Note that the two $B$ meet. 
Then we obtain the following figure.  
It is called the piece I again.

\vskip2cm
\hskip2cm\includegraphics[width=8cm]{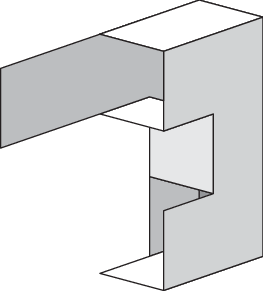}

\np

We obtain three copies of the piece I. 

Call them the first piece I, the second piece I, the third piece I.   




The $A, B, C$ is printed on each piece.

Attach the first piece I to the piece IV with the following properties.  
 $A$ meets $A$. 
 $B$ meets $B$. 
 $C$ meets $C$.  


We obtain the following.

\vskip3cm
\hskip1cm\includegraphics[width=8cm]{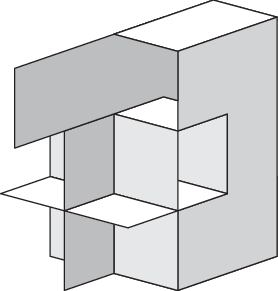}  

\np

 $A, B, C$ in the second (resp. third) piece I are called $A', B', C'$ (resp. $A'', B'', C''$).

Attach the second piece I to `the piece IV with a copy of the piece I' 
with the following properties. 
 $A$ meets $A'$.  $B$ meets $B'$.  $C$ meets $C'$.  

Attach the third piece I to `the piece IV with two copies of the piece I' 
with the following properties. 
 $A''$ meets $A'$. 
 $B''$ meets $B'$.  
 $C''$ meets $C'$.

We obtain the following figure. It is called the piece VI.
Note that the arrow in a copy of the piece I 
meets  one of the `double arrow' in another copy of the piece I.

\vskip2cm
\includegraphics[width=10cm]{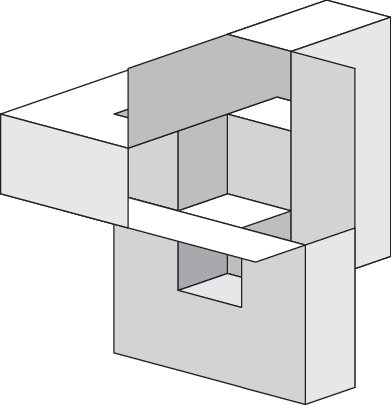}  

\np

Cut a copy of Figure II along the solid lines. What we obtain after cutting is called 
the piece II. 

Fold the piece II along the dotted line so that we see the dotted line inside, 
and make `the angle made by the paper at the dotted line' $90^\circ$.  

Use a strip of scotch tape and attach the edges which meet. 
Note that the two $P$ meet. 
Then we obtain the following.  It is called the piece II again.

\vskip1cm
\includegraphics[width=10cm]{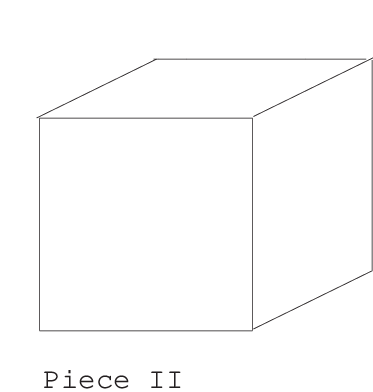}  
\vskip1cm

\np

Attach the piece II to the piece VI so that 
each star in the piece II meets each star in the piece VI. 
The result is the Boy surface. 

We draw the two figures of the Boy surface,  
seeing from two different directions.

\np

\includegraphics[width=7cm]{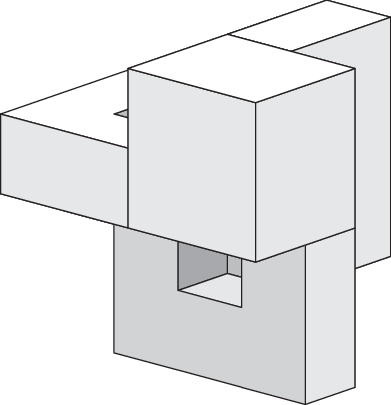}

\vskip2cm
\includegraphics[width=7cm]{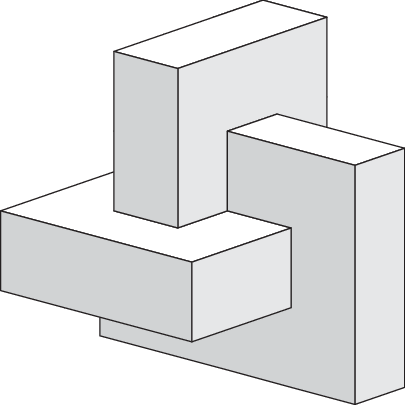}

\np

This Boy surface has the corner.  

\vskip1cm
\includegraphics[width=4cm]{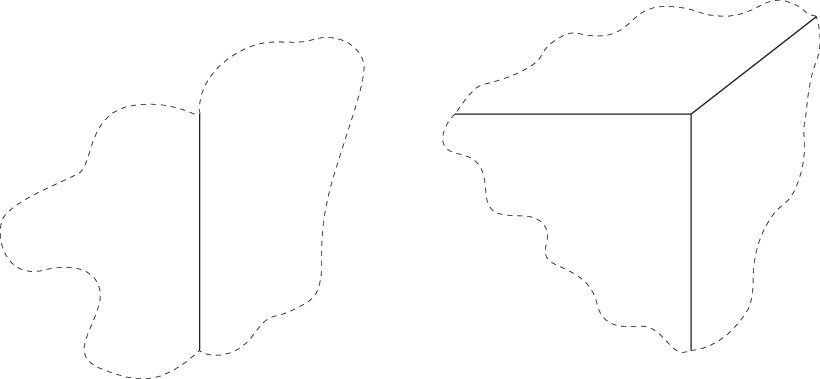}   
\vskip1cm

If you prefer the Boy surface without corner, 
imagine making the corner smooth. Or, make so.

The line which is the intersection of two sheets in the following figure  
is the set of double points.  (The double point is the intersection of two sheets.) 

\vskip1cm
\includegraphics[width=2cm]{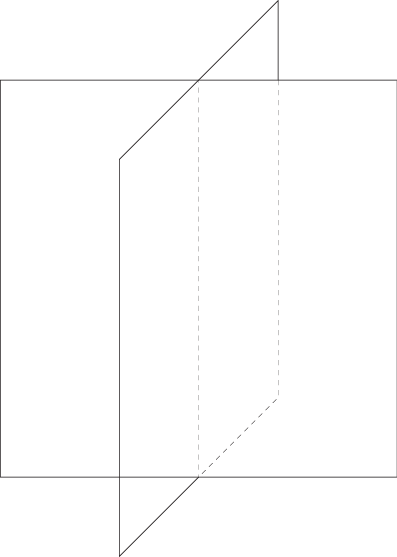} 
\vskip1cm

The point $(0,0,0)$ in the piece IV is the triple point.  
(The triple point is the intersection of three sheets as shown below.)  
The Boy surface contains only one triple point. 

\vskip1cm
\includegraphics[width=4cm]{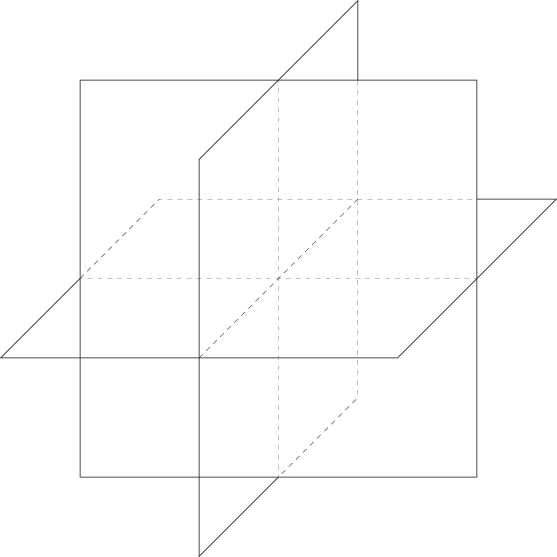}  
\vskip1cm

It is known that we cannot immerse $\R P^2$ into $\R^3$ without a triple point.

\np

\section{Prove}
Prove that the paper-craft which we made is the immersion of $\R P^2$. 

\bigbreak\noindent
Sketchy proof: 
Calculate the homology groups or the betti number of the 
manifold whose immersion is the Boy surface.  
After that, use the Poincar\'e theorem on classifying surfaces. 

\bigbreak\noindent
Alternative sketchy proof: 
Remove the piece II and the piece IV from the Boy surface. 
Then the result is made from the three copies of the piece I. 
Remove the double point set by the following procedure.  

\vskip1cm
\includegraphics[width=7cm]{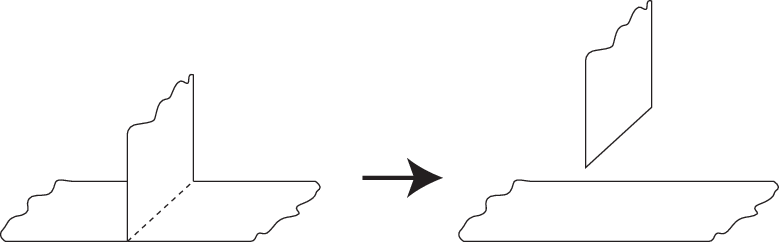}  
\vskip1cm

Prove the result is (the M\"obius band)$-$(three discs). 
Note that the boundary is a set of three circles. 
After that, use the Poincar\'e theorem on classifying surfaces.

\np

\bigbreak\noindent
Eiji Ogasa\newline 
 Computer Science, Meijigakuin University, Yokohama, Kanagawa, 244-8539, Japan 
\newline 
ogasa@mail1.meijigkakuin.ac.jp 
\quad
pqr100pqr100@yahoo.co.jp  

\bb

The author made a movie to explain our paper-craft. He put it in Youtube. 
The website of the movie is connected with the author's website. 
One can find these websites by typing in 
the author's name, `Eiji Ogasa', 
or 
the title of this article, `Make your Boy surface',   
in the search engine.


\begin{thebibliography}{ABCD}
\bibitem{Oni}
E. Ogasa: 
Boy surface sekkeisho (A translation of this article by the author. 
One can find it by typing in 
the author's name, `Eiji Ogasa', or the Japanese title.)
\end{thebibliography}
\end{document}